\documentclass[a4paper,12pt]{article}
\usepackage[utf8]{inputenc}

\usepackage{amsmath}
\usepackage{amssymb,latexsym}
\usepackage{amsthm}

\usepackage[top=3cm, left=3cm, right=2cm, bottom=3cm]{geometry}
\newtheorem{theorem}{Theorem}
\newtheorem{lemma}[theorem]{Lemma}
\newtheorem{proposition}[theorem]{Proposition}

\title{Explicit bounds for Lipschitz constant of solution to basic problem in calculus of variations}
\author{Miguel Oliveira, Georgi Smirnov\\~\\
University of Minho}
\date{~}

\begin{document}

\maketitle

\begin{abstract}
In this paper we present explicit estimate for Lipschitz constant of solution to a problem of calculus of variations.
The approach we use is due to Gamkrelidze and is based on the equivalence of the  problem of calculus of variations and a time-optimal control problem. The obtained estimate is used to compute  complexity bounds for a path-following method applied to a convex problem of calculus of variations with polyhedral end-point constraints.
\end{abstract}


\section{Introduction}

The first works concerning regularity of solutions to basic problem in calculus of variations appeared more than a century ago  \cite{B,Tonelli}.
In the last 30 years, regularity of solution to problems of calculus of variations and optimal control has been a subject of intensive studies  
(see, e.g.,  \cite{ClV,Cellina,CF,CFM,Z}). However, at least to our knowledge, there are no explicit bounds for Lipschitz constant of solution to the problem of calculus of variations.
 In this paper, we obtain such a bound. This is done under rather strong conditions. The approach we use is very close to the one from \cite{ST,T} and is based on the equivalence between the problem of calculus of variations  and the time optimal control problem established by Gamkrelidze \cite{G}.
 In the case of basic problem of calculus of variations, for lagrangians considered here, the lipschitzian regularity  was proved by Clarke and Vinter in \cite{ClV} but  their proof  is not constructive: they do not obtain an explicit expression for the Lipschitz constant. 
 
 Based on our explicit estimate for the Lipschitz constant of solution we derive complexity bounds for a path-following method applied to a convex problem of calculus of variations with polyhedral end-point constraints. It is well-known from the approximation theory that the question about how functions can best be approximated with simpler functions is closely related with their regularity \cite{KT}. We  use the regularity properties  as main tool to approximate convex  problems of calculus of variations by a convex programming problem and to get the respective complexity bounds for a path-following method \cite{NN}.

\section{Main results}

Let us introduce some notations used in the sequel. 
We denote the norm of the vector $x\in R^n$ by $|x|$ and the inner product of two vectors $x_1,x_2\in R^n$ by $\langle x_1,x_2\rangle$. The closed unit ball in $R^n$ is denoted by $B_n$.  The distance between $x\in R^n$ and $C\subset R^n$ is denoted by $d(x,C)$. The convex hull of $C$ is denoted by co$C$. The tangent cone to $C\subset R^n$ at $x\in C$ is defined as ${\cal T}(C,x)=\{ v\in R^n\mid \lim_{\lambda\downarrow 0}\lambda^{-1}d(x+\lambda v,C)=0\}$. The conjugate cone to a cone $K\subset R^n$ is denoted by $K^*=\{ x^*\in R^n \mid \langle x^*,x\rangle\geq 0,\; x\in K\}$. The set of absolutely continuous functions $x:[t_0,t_1]\rightarrow R^n$ is denoted by $AC([t_0,t_1],R^n)$ and the set of measurable essentially bounded functions $x:[t_0,t_1]\rightarrow R^n$ is denoted by $L_{\infty}([t_0,t_1],R^n)$. The graph of a set-valued map $F:R^m\rightarrow R^n$ is denoted by gr$F$.

\subsection{Lipschitzian regularity of solutions}

Consider the following problem of calculus of variations
\begin{eqnarray}
\label{1}
&& \int_0^1 L(t,x(t),\dot{x}(t))dt \rightarrow\inf,\\
&& x(0)=0,\;\;\; Ax(1)\leq b, \nonumber
\end{eqnarray}
where $L:R\times R^n\times R^n\rightarrow R$, $A$ is an $m\times n$-matrix, and $b\in R^m$. We shall denote the set $\{ x\mid Ax\leq b\}$ by $S$.

Assume that the following conditions are satisfied:

\vspace{5mm}

\noindent (C1) The function $L(\cdot,\cdot,\cdot)$ is  continuously differentiable and satisfies the coercivity condition
\begin{equation}
\label{gr1}
L(t,x,u)\geq \theta (|u|)>0,
\end{equation}
where $\lim_{r\rightarrow\infty}{r}/{\theta(r)}=0$.

\noindent (C2) The function $L(t,x,\cdot)$ is strictly convex, i.e., there exists a constant $\mu>0$ such that
$$
L(t,x,u)+\langle\nabla_u L(t,x,u),v-u\rangle +\frac{\mu}{2}|v-u|^2\leq L(t,x,v)
\;\;{\rm for\: all}\;\; t,x,u,v.
$$
\noindent (C3) There exist constants $\xi>0$ and $\delta>0$ such that
$$
|\nabla_{(t,x)}L(t,x,u)|\}\leq \xi L(t,x,u)+\delta.
$$

\vspace{5mm}

Let $r_0>0$ be such that $\theta(r)/r\geq 1$, whenever $r\geq r_0$, and let $a\in S$.
Set
$$
c=r_0+\int_0^1L(t,ta,a)dt.
$$
Let $\hat{x}(\cdot)\in AC([0,1],R^n)$ be the solution to problem of calculus of variations (\ref{1}).
 Denote by $M$ the set of points $t\in[0,1]$ such that $|\dot{\hat{x}}(t)|\leq \theta(|\dot{\hat{x}}(t)|)$.
Since
$$
 c\geq r_0+ \int_0^1 L(t,\hat{x}(t),\dot{\hat{x}}(t))dt\geq r_0+ \int_M L(t,\hat{x}(t),\dot{\hat{x}}(t))dt
$$
\begin{equation}
 \geq \int_{[0,1]\setminus M} r_0dt+\int_M\theta(|\dot{\hat{x}}(t)|)dt\geq \int_0^1|\dot{\hat{x}}(t)|dt,
\label{ner1}
\end{equation}
we have
$$
|\hat{x}(t)|\leq c,
$$ 
whenever $t\in [0,1]$. Set $\Omega=[0,1]\times c B_n$.

Put 
$$
\Lambda_0=\max_{(t,x)\in\Omega}L(t,x,0),
$$ 
$$
\Lambda_1=\max_{(t,x)\in\Omega}|\nabla_u L(t,x,0)|,
$$ 
and 
$$
\sigma (r) = \max_{(t,x)\in\Omega,\; |u|\leq r} (\langle\nabla_uL(t,x,u),u\rangle-L(t,x,u)).
$$ 
The following Lemma is an immediate consequence of (C2).

\begin{lemma}
\label{lem00}
The following inequality holds:
\begin{equation}
\label{La}
-\Lambda_0-\Lambda_1|u|+\frac{\mu}{2}|u|^2\leq L(t,x,u),\;\;\; (t,x)\in \Omega.
\end{equation}
The function $\sigma (r)$ tends to infinity as $r\rightarrow\infty$.
\end{lemma}

{\it Proof}. Inequality (\ref{La}) is a consequence of the following inequality:
$$
L(t,x,0)+\langle\nabla_u L(t,x,0),u\rangle + \frac{\mu}{2}|u|^2\leq L(t,x,u).
$$
Next, from condition (C2) we have
$$
L(t,x,u)-\langle\nabla_u L(t,x,u),u\rangle+\frac{\mu}{2}|u|^2\leq L(t,x,0).
$$
Therefore
$$
\frac{\mu}{2}r^2-\Lambda_0\leq \sigma (r).\;\;\; \qed
$$

\vspace{5mm}

By Lemma \ref{lem00} there exists $0<T_0<1$ such that $\beta=\sigma((c+1)/T_0)>\delta/\xi$.
Put
$$
\eta=\sup_{r\geq 0}\frac{r}{\theta(r)+\beta},
$$
\begin{equation}
\label{gam}
\bar{\gamma}=e^{\eta\xi(c+T_0\beta)},
\end{equation}
and
$$
\Lambda_2=\max_{(t,x)\in\Omega,\: |u|=\bar{\gamma}+1}L(t,x,u)
$$ 

 Let $\varrho>0$ be such that
\begin{equation}
\label{rh}
\frac{r}{\left(\frac{\mu}{2}r^2-\Lambda_1r-\Lambda_0\right) +\beta}<\frac{1}{\Lambda_2+\beta},
\end{equation}
whenever $r\geq\varrho$.

\begin{theorem}\label{th1}
The following inequality holds:
$$
|\dot{\hat{x}}(t)|\leq \ell=\max\left\{ \varrho,\sqrt{\frac{2}{\mu}\left(\Lambda_0+\beta\right)},\frac{\Lambda_1+\sqrt{\Lambda_1^2+4\mu\Lambda_0}}{2}\right\}.
$$
\end{theorem}

\subsection{Complexity bounds}

Recall the path-following method from \cite[Ch. 7]{NN}. Let $P:R^n\rightarrow R$  be a convex function. Consider the problem
\begin{equation}\label{PP}
\left. \begin{array}{c} 
\sigma \rightarrow \min, \\
x=(x^{(1)}, x^{(2)}, \ldots, x^{(n)}) \in R^n, \\
P(x)\leq \sigma, \\
Ax\leq b, \\
 |x^{(i)}|\leq M, \ \ \ i=\overline{1,n}.
\end{array} \right\}
\end{equation}
Let $F$ be the function defined by
$$
F(x, \sigma)=-\ln(\sigma-P(x))-\ln(\bar{\sigma}-\sigma)
$$
$$
-\sum_{i=1}^m\ln(b_i-(Ax)^{(i)})-\sum_{i=1}^n\ln(M^2-(x^{(i)})^2),
$$
where $\bar{\sigma}\geq\max_{\{x||x^{(i)}|\leq M, i=\overline{1,n}\}}P(x)$.
We use the notation 
$$
\| v\|_{x}^F=\langle [\nabla^2F(x)]^{-1}v,v \rangle^{1/2}.
$$
Let $\kappa \in (0,1)$ and let $\gamma>0$ be such that $\gamma\leq \frac{\sqrt{\kappa}}{1+\sqrt{\kappa}}-\kappa$.

\

\emph{Path-following Method}
\begin{itemize}
\item Initialization: Set $\alpha_0=0$ and $v=(0,1)\in R^{n+1}$. Choose an accuracy $\varepsilon>0$, $x_0\in R^n$, and $\sigma_0\in R$ such that
$$
\| \nabla F(x_0,\sigma_0)\|_{(x_0,\sigma_0)}^F \leq \kappa.
$$

\item Step $k$: Set
$$
\alpha_{k+1}=\alpha_k+\frac{\gamma}{\| v\|_{(x_k,\sigma_k)}^F},
$$
$$
(x_{k+1},\sigma_{k+1})=(x_k,\sigma_k)- [\nabla^2F(x_k,\sigma_k)]^{-1}(\alpha_{k+1}v+ \nabla F(x_k,\sigma_k)).
$$

\item Stop the process if 
$$
m+n+1+\frac{(\kappa+\sqrt{m+n+1})\kappa}{1-\kappa}\leq \varepsilon \alpha_k.
$$

\end{itemize}

Let $\mathcal{N}$ be the largest integer satisfying
$$
\mathcal{N}\leq \frac{\ln\left( \frac{(1+\kappa)(m+n+1)+(\kappa+\sqrt{m+n+1})\kappa}{\gamma(1-2\kappa)\varepsilon}\|v\|_{(\tilde x, \tilde \sigma)}^F\right)}{\ln\left (1+\frac{\gamma}{\kappa+\sqrt{m+n+1}}\right)}+1 ,
$$
where $(\tilde x,\tilde \sigma)={\rm argmin}\: F$.
Since $\|v\|_{(\tilde x, \tilde \sigma)}^F$ is the $(n+1,n+1)$-th entry of matrix $[\nabla^2F(\tilde{x}, \tilde \sigma)]^{-1}$, using the Sherman-Morrison-Woodbury formula, we obtain
$$
\|v\|_{(\tilde{x}, \tilde \sigma)}^{F}\leq \left(\frac{1}{(\tilde{\sigma}-P(\tilde{x}))^2}+\frac{1}{(\bar{\sigma}-\tilde{\sigma})^2}\right)^{-1/2}
$$
$$
=\left(\frac{(\tilde{\sigma}-P(\tilde{x}))^2(\bar{\sigma}-\tilde{\sigma})^2}{(\tilde{\sigma}-P(\tilde{x}))^2+(\bar{\sigma}-\tilde{\sigma})^2}\right)^{1/2}.
$$
Let us consider the function
$$
g(\lambda)=\frac{(\lambda-P(\tilde{x}))^2(\bar{\sigma}-\lambda)^2}{(\lambda-P(\tilde{x}))^2+(\bar{\sigma}-\lambda)^2}, \ \ \lambda\in \left[P(\tilde{x}), \bar{\sigma}\right].
$$
Its maximum is achieved at $\lambda=(\bar{\sigma}+P(\tilde{x}))/2$. Thus, we have
$$
\|v\|_{(\tilde{x}, \tilde \sigma)}^F\leq\frac{1}{2\sqrt{2}}\left(\bar{\sigma}+P(\tilde{x})\right)
$$
$$
\leq \frac{1}{2\sqrt{2}}\left(\bar{\sigma}+\max_{\{x||x^{(i)}|\leq M, i=\overline{1,n}\}}P(x)\right) \leq\frac{1}{\sqrt{2}}\bar{\sigma}.
$$

\begin{proposition}{\cite[Ch. 7]{NN}}
\label{proposition5}
The path-following method terminates no more than after $\mathcal{N}$ steps.  
At the moment of termination we have $|P(x_\mathcal{N})-P(\hat x)|< \varepsilon$, 
where $\hat x$ is the solution of problem (\ref{PP}).
\end{proposition}

Assume that $L$ is strictly convex function of $(x,u)$:
$$
L(t,x_1,u_1)+\langle \nabla_x L(t,x_1,u_1),x_2-x_1\rangle+\langle \nabla_u L(t,x_1,u_1),u_2-u_1\rangle
$$
\begin{equation}
\label{sc}
+\frac{\mu}{2}(|x_2-x_1|^2+|u_2-u_1|^2)\leq L(t,x_2,u_2),\;\;{\rm for\: all}\;\; t,x_1,x_2,u_1,u_2.
\end{equation}
Put
$$
K_L=\max_{t\in [0,1], |x|\leq\ell,  |u|\leq\ell} |\nabla_{(t,x,u)}L(t,,x,u)|.
$$
We denote by $\mathcal{U}_N\subset L_\infty([0,1],R^n)$ the set  consisting of piecewise constant functions $u(\cdot)$ taking values $u(t) = u(\tau k)$, $t\in ]\tau k, \tau (k + 1)]$, $k = \overline{0,N-1}$, $\tau = 1/N$ and satisfying $|u(t)|\leq \ell$.

Let $\mathcal{F}:R^{n\times N}\times R\rightarrow R$ be the function defined by
$$
\mathcal{F}(u,\sigma)=-\ln\left(\sigma-
\int_0^1 L\left(t,\int_0^tu(s)ds,u(t)\right)dt
\right)-\ln(\bar{\sigma}-\sigma)
$$
$$
-\sum_{i=1}^m\ln\left(b_i-\left(A\sum_{k=0}^{N-1}u(\tau k) \right)^{(i)}\right)-\sum_{i=1}^n\sum_{k=0}^{N-1}\ln(\ell^2-(u^{(i)}(\tau k))^2),
$$
where $u(\cdot)\in \mathcal{U}_N$ and 
\begin{equation}
\label{sig}
\bar{\sigma}=L(0,0,0)+K_L(1+2\ell)
\geq\max_{\{u||u|\leq \ell\}} \int_0^1 L\left(t,\int_0^t u(s)ds,u(t)\right)dt.
\end{equation}
\begin{theorem}
\label{th3}
Let $\varepsilon>0$,
$$
N>\frac{4(1+\ell ) K_L}{\varepsilon}, 
$$
and
$$
\mathcal{N}\geq \frac{\ln\left(\sqrt{2}\frac{(1+\kappa)(m+nN+1)+(\kappa+\sqrt{m+nN+1})\kappa}{\gamma(1-2\kappa)\varepsilon}\bar{\sigma}\right)}{\ln\left(1+\frac{\gamma}{\kappa+\sqrt{m+nN+1}}\right)}+1.
$$
The path-following method with the function $\cal F$ finds an admissible  $\breve{u}\in \mathcal{U}_N$ and the respective trajectory 
$$
\breve{x}(t)=\int_0^t\breve{u}(s)ds
$$
satisfying 
$$
\left|\int_0^1 (L(t,\hat{x}(t), \dot{\hat{x}}(t))-L\left(t,\breve{x}(t),\breve{u}(t))\right)dt\right|< \varepsilon
$$ 
and 
$$
\int_0^1 \left(\left|\breve{x}(t)-\hat{x}(t)\right|^2+|\breve{u}(t)-\dot{\hat{x}}(t)|^2\right)dt\leq
\frac{2}{\mu}
\varepsilon,
$$ no more than after $\mathcal{N}$ iterations.
\end{theorem}

\section{Proofs of main results}

First, recall some results needed to prove Theorem \ref{th1}.

\subsection{Background notes}

According to Gamkrelidze \cite{G} problem (\ref{1}) is equivalent to the time-optimal control problem
\begin{eqnarray}
&& T\rightarrow\inf,\nonumber\\
&& \frac{d(t,y)(\tau )}{d\tau}
=
\frac{(1,w(\tau ))}{L(t(\tau ),y(\tau ),w(\tau ))},\label{2}\\
&& (t,y)(0)=(0,0),\;\;\; (t,y)(T )\in (1,S)\nonumber,
\end{eqnarray}
where
$$
\tau(t)=\int_0^t L(s,x(s),\dot{x}(s))ds,
$$
 $y(\tau)=x(t(\tau))$, and $w(\tau)=\dot{x}(t(\tau))$. Note that, since the function $\tau=\tau(t)$ is strictly monotonous and absolutely continuous, its inverse, $t=t(\tau)$, is
also strictly monotonous and absolutely continuous. Therefore $y=y(\tau)$ is absolutely continuous and $w=w(\tau)$ is measurable. 
Recall the following proposition from \cite[Sec. 8.5]{G}.

\begin{proposition}
\label{pr1}
Assume that $L$ is continuously differentiable and   condition (\ref{gr1}) is satisfied. Then
\begin{enumerate}
\item For any function $x(\cdot )\in AC([t_0,t_1],R^n)$ there exists a trajectory $(t,y)(\tau )$, $\tau\in [0,T]$ of control system (\ref{2}) such that
$(t,y)(0)=(t_1,x(t_1))$, $(t,y)(T)=(t_2,x(t_2))$, and
$$
T=\int_{t_1}^{t_2} L(t,x(t),\dot{x}(t))dt.
$$
\item For any trajectory $(t,y)(\tau )$, $\tau\in [0,T]$ of control system (\ref{2}) such that $\frac{d}{d\tau}(t,y)\neq (0,0)$ almost everywhere, there exists 
$x(\cdot )\in AC([t(0),t(T)], R^n)$ such that $x(t(0))=y(0)$, $x(t(T))=y(T)$, and
$$
T=\int^{t(T)}_{t(0)} L(t,x(t),\dot{x}(t))dt.
$$
\end{enumerate}
\end{proposition}

Let $F:R^n\rightarrow R^n$ be a Lipschitzian set-valued map with compact values. Consider the following time-optimal problem
\begin{eqnarray*}
&& T\rightarrow\min,\\
&& \dot{x}\in F(x),\\
&& x(0)=x_0,\;\; x(T)\in S.
\end{eqnarray*}
Here $S\subset R^n$ is a closed convex set. 
Let $\hat{x}(\cdot)\in AC([0,T],R^n)$ be a solution to this problem. Consider a convex cone $K(t)\subset {\cal T}({\rm gr\:co}F,(\hat{x}(t),\dot{\hat{x}}(t)))$ measurably depending on $t\in [0,T]$.

There exist many necessary conditions of optimality for time-optimal problems with differential inclusions (see e.g. \cite{Cl,V}).
For our considerations the most suitable formulation is contained in the  following proposition which is a consequence of  \cite[Theorem 5]{S1}.
\begin{proposition}
\label{pr2}
There exists a function $p(\cdot )\in AC([0,T],R^n)$ such that
\begin{enumerate}
\item $(\dot{p},p)\in -K^*(t)$, $\langle p(t),\dot{\hat{x}}(t)\rangle\equiv h\geq 0$;
\item $p(T)\in ({\cal T}(S,\hat{x}(T)))^*$;
\item $|p(T)|>0$.
\end{enumerate}
\end{proposition}

In the case of a smooth control system, the Yorke approximation can be chosen as the cone $K(t)$.
Let $U\subset R^k$, and let $f:R^n\times U\rightarrow R^n$ be a function. Assume that $f$ is differentiable in $x$ and the set $f(x,U)$ is convex for all $x\in R^n$. For $(\hat{x},\hat{u})\in R^n\times U$ denote $\hat{v}=f(\hat{x},\hat{u})$ and set $C=\nabla_x f(\hat{x},\hat{u})$ and $K={\cal T}(f(\hat{x},U),\hat{v})$. Recall the following proposition \cite[p. 38]{S2}.
\begin{proposition}
\label{pr3}
The following inclusion holds:
$$
\{ (x,v)\in R^n\times R^n\mid v\in Cx+K\}\subset {\cal T}((\hat{x},\hat{u}), {\rm gr} f(\cdot,U)).
$$
\end{proposition}

Recall also the following useful formula \cite[p. 50]{S2}.
\begin{proposition}
\label{pr4}
Let $C:R^n\rightarrow R^n$ be a linear operator, and let $K\subset R^n$ be a convex cone. Then the following equality holds:
$$
\{ (x,v)\in R^n\times R^n\mid v\in Cx+K\}^*
$$
$$
=
\{ (x^*,v^*)\in R^n\times R^n\mid x^*=-C^*v^*,\; v^*\in K^*\}.
$$
\end{proposition}

\subsection{Proof of Theorem \ref{th1}}

Let  us consider the set-valued map
$$
G(t,y)=\left\{ (v^0,v)\in R\times R^n\mid 
v^0=\frac{\rho }{ L(t,y,q)+\beta},\; \right.
$$
$$
\left. v=\frac{\rho w}{ L(t,y,w)+\beta},\; w\in U,\; \rho\in [0,1]\right\}.
$$
\begin{lemma}\label{lem1}
The set-valued map $G$ has convex compact values and is Lipschitzian  in the set $\Omega$.
\end{lemma}

{\it Proof}. Let  $(v^0_i,v_i)\in G(t,y)$, $i=1,2$. There exist $\rho_i\in [0,1]$ and $w_i\in U$, $i=1,2$, such that
$$
v_i^0=\frac{\rho_i}{ L(t,y,w_i)+\beta},\;\;\;
 v_i=\frac{\rho_i w_i}{ L(t,y,w_i)+\beta},\;\;\; i=1,2.
$$
Let  $\alpha_1,\; \alpha_2 \geq 0$, $\alpha_1+\alpha_2 =1$. Show that 
 $\alpha_1 (v^0_1,v_1) +\alpha_2 (v^0_2,v_2)\in G(t,y)$. Put 
$$
\alpha'_1= \frac{\alpha_1\rho_1}{ L(t,y,w_1)+\beta}
\left( \frac{\alpha_1\rho_1}{ L(t,y,w_1)+\beta}
+\frac{\alpha_2\rho_2}{ L(t,y,w_2)+\beta}
\right)^{-1} 
$$
 and
$$
\alpha_2'= \frac{\alpha_2\rho_2}{ L(t,y,w_2)+\beta}
\left( \frac{\alpha_1\rho_1}{ L(t,y,w_1)+\beta}
+\frac{\alpha_2\rho_2}{L (t,y,w_2)+\beta}
\right)^{-1}. 
$$
Obviously $\alpha'_1,\; \alpha_2' \geq 0$, $\alpha'_1+\alpha_2' =1$. 
Set $w=\alpha_1'w_1+\alpha_2'w_2$. Then we get $ L(t,y,w)\leq\alpha_1'L(t,y,w_1)+\alpha_2' L(t,y,w_2)$.
We have
$$
\alpha_1 v_1+\alpha_2 v_2=
\frac{\alpha_1\rho_1 w_1}{ L(t,y,w_1)+\beta}+
\frac{\alpha_2\rho_2 w_2}{ L(t,y,w_2)+\beta}
$$
$$
=
\left(\frac{\alpha_1\rho_1}{ L(t,y,w_1)+\beta}+
\frac{\alpha_2\rho_2}{ L(t,y,w_2)+\beta}\right) w.
$$
Put
$$
\rho=
\left(\frac{\alpha_1\rho_1}{ L(t,y,w_1)+\beta}+
\frac{\alpha_2\rho_2}{ L(t,y,w_2)+\beta}\right)( L(t,y,w)+\beta).
$$
Then we obtain
$$
\rho\leq 
\frac{\alpha_1\rho_1}{ L(t,y,w_1)+\beta}( L(t,y,w_1)+\beta)+
\frac{\alpha_2\rho_2}{ L(t,y,w_2)+\beta}(  L(t,y,w_2)+\beta)
$$
$$
=\alpha_1\rho_1+\alpha_2\rho_2\leq 1,
$$
i.e., $\rho\in [0,1]$. Therefore $G(t,y)$ is convex.

From (\ref{gr1}) we have
\begin{equation}
\label{int_***}
|v^0|\leq
\frac{1}{ \theta (|w|)+\beta},\;\;\;
|v|\leq
\frac{|w|}{ \theta (|w|)+\beta},\;\;\;
\forall\: (v^0,v)\in G(t,y).
\end{equation}
Let
$$
(v_k^0,v_k)=
\left(\frac{\rho_k }{ L(t,y,w_k)+\beta},\frac{\rho_k w_k}{ L(t,y,w_k)+\beta}\right),
$$
where $w_k\in R^n$ and $\rho_k\in [0,1]$, be a sequence converging to a point  
$(v^0_0,v_0)$. If the sequence $w_k$ is bounded, then, without loss of generality, the sequence $(w_k,\rho_k)$ converges. Passing to the limit we get
 $(v^0_0,v_0)\in G(t,y)$. If the sequence $w_k$ is unbounded, then there exists a subsequence converging to infinity. Without loss of generality  $w_k$ goes to infinity. From inequalities   (\ref{int_***}) we obtain $(w_k^0,w_k)\rightarrow (0,0)$. Hence  $(w_0^0,w_0)=(0,0)\in G(t,y)$. Thus   $G(t,y)$ is a closed set. From (\ref{int_***}) we see that it is bounded.

Let $(t_1,y_1)$ and $(t_2,y_2)$ be two points in $\Omega$. 
Let 
$$
 (v_1^0,v_1)=\frac{(\rho,\rho w)}{ L(t_1,y_1,w)+\beta}\in G(t_1,y_1).
$$
Consider the point
$$
 (v_2^0,v_2)=\frac{(\rho,\rho w)}{ L(t_2,y_2,w)+\beta}\in G(t_2,y_2).
$$
Since $\beta>\delta/\xi$, from (C3) we have
$$
|v^0_1-v_2^0|\leq\max_{\lambda\in [0,1]}\frac{ |\nabla_{(t,x)}L(\lambda t_1+(1-\lambda )t_2,\lambda y_1+(1-\lambda )y_2,w)|}{( L(\lambda t_1+(1-\lambda )t_2,\lambda y_1+(1-\lambda )y_2,w)+\beta)^2}(|t_1-t_2|+|y_1-y_2|)
$$
$$
\leq\frac{\xi}{\beta} (|t_1-t_2|+|y_1-y_2|)
$$
and
$$
|v_1-v_2|\leq\max_{\lambda\in [0,1]}\frac{ |\nabla_{(t,x)}L(\lambda t_1+(1-\lambda )t_2,\lambda y_1+(1-\lambda )y_2,w)||w|}{( L(\lambda t_1+(1-\lambda )t_2,\lambda y_1+(1-\lambda )y_2,w)+\beta)^2}(|t_1-t_2|+|y_1-y_2|)
$$
$$
\leq\eta\xi (|t_1-t_2|+|y_1-y_2|),
$$
i.e. $G$  is Lipschitzian  in the set $\Omega$. \qed

\vspace{5mm}

Let $\hat{x}(\cdot)$ be a solution to problem (\ref{1}). By the first part of Proposition \ref{pr1} there exists
a trajectory $(\hat{t},\hat{y})(\tau )$, $\tau\in [0,\hat{T}]$, of control system 
\begin{equation}\label{3}
\frac{d(t,y)(\tau )}{d\tau}
=
\frac{(1,w(\tau ))}{ L(t(\tau ),y(\tau ),w(\tau ))+\beta},\;\; w(\tau)\in R^n.
\end{equation}
 such that
$(\hat{t},\hat{y})(0)=(0,0)$, $(\hat{t},\hat{y})(\hat{T})=(1,a)$, and
$$
\hat{T}=\int_{0}^{1} ( L(t,\hat{x}(t),\dot{\hat{x}}(t))+\beta)dt.
$$
The control corresponding to $\hat{y}(\cdot)$ is denoted by $\hat{w}(\cdot)$. 

\begin{lemma}
\label{lem10}
There exists a nonzero function  $(q,p)(\cdot)\in AC([0,\hat{T}],R\times R^n)$ such that
\begin{eqnarray}
&& \frac{dq}{d\tau}=\frac{ (q+\langle \hat{w},p\rangle )L_t(\hat{t}(\tau ),\hat{y}(\tau ),\hat{w}(\tau )) }{( L(\hat{t}(\tau ),\hat{y}(\tau ),\hat{w}(\tau ))+\beta)^2},\label{P1}\\
&& \frac{dp}{d\tau}= \frac{ (q+\langle \hat{w},p\rangle )\nabla_x L(\hat{t}(\tau ),\hat{y}(\tau ),\hat{w}(\tau )) }{( L(\hat{t}(\tau ),\hat{y}(\tau ),\hat{w}(\tau ))+\beta)^2},\label{P2}\\
&&\label{P3}
\frac{p}{ L(\hat{t}(\tau ),\hat{y}(\tau ),\hat{w}(\tau ))+\beta}
-\frac{(q+\langle \hat{w},p\rangle )\nabla_w L(\hat{t}(\tau ),\hat{y}(\tau ),\hat{w}(\tau ))}{( L(\hat{t}(\tau ),\hat{y}(\tau ),\hat{w}(\tau ))+\beta)^2}=0,\\
&&\label{P4} \frac{ q+\langle \hat{w},p\rangle  }{ L(\hat{t}(\tau ),\hat{y}(\tau ),\hat{w}(\tau ))+\beta}\equiv h>0.
\end{eqnarray}
\end{lemma}

{\it Proof}. From the second part of Proposition \ref{pr1} we see that  $(\hat{t},\hat{y})(\tau )$, $\tau\in [0,\hat{\tau}]$ is a solution to the problem 
\begin{eqnarray}
&& T\rightarrow\inf,\nonumber\\
&& \frac{d(t,y)(\tau )}{d\tau}
=
\frac{(1,w(\tau ))}{ L(t(\tau ),y(\tau ),w(\tau ))+\beta},\;\; w(\tau)\in R^n, \label{4}\\
&& (t,y)(0)=(0,0),\;\;\; (t,y)(T )=(1,a).\nonumber
\end{eqnarray}
The time-optimal problem
\begin{eqnarray}
&& T\rightarrow\inf,\nonumber\\
&& \frac{d(t,y)(\tau )}{d\tau}\in G(t,y),\label{TG}\\
&& (t,y)(0)=(0,0),\;\;\; (t,y)(T )=(1,a)\nonumber,
\end{eqnarray}
also has a solution $(\tilde{t},\tilde{y})(\tau)$, $\tau\in [0,\tilde{T}]$. By the Filippov lemma there exists a measurable function $(\tilde{\rho},\tilde{w}) (\tau)$, $\tau\in [0,\tilde{T}]$, such that
$$
\frac{d(\tilde{t},\tilde{y})(\tau )}{d\tau}
=
\frac{\tilde{\rho}(\tau)(1,\tilde{w}(\tau ))}{ L(\tilde{t}(\tau ),\tilde{y}(\tau ),\tilde{w}(\tau ))+\beta}
$$
at almost all points where $d(\tilde{t},\tilde{y})/d\tau\neq (0,0)$. 
Applying Propositions \ref{pr2}-\ref{pr4}, we see that there exist  $(q,p)(\cdot)\in AC([0,\tilde{T}],R\times R^n)$, a nonzero function, and a constant $h\geq 0$ such that
\begin{eqnarray}
&& \frac{dq}{d\tau}=\frac{ \tilde{\rho}(q+\langle \tilde{w},p\rangle )L_t(\tilde{t}(\tau ),\tilde{y}(\tau ),\tilde{w}(\tau )) }{( L(\tilde{t}(\tau ),\tilde{y}(\tau ),\tilde{w}(\tau ))+\beta)^2},\label{n1}\\
&& \frac{dp}{d\tau}=\frac{ \tilde{\rho}(q+\langle \tilde{w},p\rangle )\nabla_x L(\tilde{t}(\tau ),\tilde{y}(\tau ),\tilde{w}(\tau )) }{( L(\tilde{t}(\tau ),\tilde{y}(\tau ),\tilde{w}(\tau ))+\beta)^2},\label{n2}\\
&& h\equiv\frac{\tilde{\rho}(q+\langle \tilde{w},p\rangle )}{ L(\tilde{t}(\tau ),\tilde{y}(\tau ),\tilde{w}(\tau ))+\beta}\nonumber\\
&&
\geq
\frac{{\rho}(q+\langle {w},p\rangle )}{ L(\tilde{t}(\tau ),\tilde{y}(\tau ),{w})+\beta},\;\;\rho\in [0,1],\;\; w\in R^n,\label{n3}
\end{eqnarray}
at almost all points such that $d(\tilde{t},\tilde{y})/d\tau\neq (0,0)$. If $d(\tilde{t},\tilde{y})/d\tau = (0,0)$ on a set of positive measure, then $h=0$. At points where $\tilde{\rho} (\tau)>0$, from maximum condition (\ref{n3}) we have
\begin{equation}
\label{n31}
\frac{\tilde{\rho}p}{ L(\tilde{t}(\tau ),\tilde{y}(\tau ),\tilde{w}(\tau ))+\beta}
-\frac{ \tilde{\rho}(q+\langle \tilde{w},p\rangle )\nabla_w L(\tilde{t}(\tau ),\tilde{y}(\tau ),\tilde{w}(\tau ))}{( L(\tilde{t}(\tau ),\tilde{y}(\tau ),\tilde{w}(\tau ))+\beta)^2}=0.
\end{equation}
Since $q+\langle \tilde{w},p\rangle =0$, we get $p=0$ and $q=0$, a contradiction. Thus $\tilde{\rho}(\tau)=0$ almost everywhere. This is impossible. Hence $d(\tilde{t},\tilde{y})/d\tau\neq (0,0)$ at almost all points $\tau\in [0,\tilde{\tau}]$. Therefore conditions (\ref{n1})-(\ref{n3}) are satisfied almost everywhere and $\tilde{\rho}(\tau)>0$ at almost all points $\tau\in [0,\tilde{\tau}]$. Thus $h>0$, because the equality $h=0$ implies, as above, $(q,p)(\tau)\equiv 0$. From   (\ref{n3}) we obtain $\tilde{\rho}(\tau)=1$. Thus we can identify the trajectories $(\hat{t},\hat{y})(\cdot)$ and $(\tilde{t},\tilde{y})(\cdot)$. Both of them are solutions to time-optimal problem (\ref{TG}) and satisfy necessary conditions of optimality (\ref{n1}), (\ref{n2}), and (\ref{n31}) with $\tilde{\rho}=1$. \qed

\vspace{5mm}
Denote by $\hat{\tau}(\cdot)$ the function inverse to $\hat{t}(\cdot)$. Then we have 
$\dot{\hat{x}}(\cdot)=\hat{w}(\hat{\tau}(\cdot))$. Therefore it suffices to obtain the bounds for $\hat{w}(\cdot)$.
 We shall use the notation $\hat{L}(\tau)$ for $L(\hat{t}(\tau),\hat{y}(\tau),\hat{w}(\tau))$.

\begin{lemma}
\label{lem2}
If $q(\tau)\leq 0$, then $|\hat{w}(\tau)|>(c+1)/T_0$.
\end{lemma}

{\it Proof}. Multiplying (\ref{P3}) by $\hat{w}(\tau)$, we obtain
\begin{equation}
\label{P44}
( \hat{L}(\tau)+\beta)\langle \hat{w}(\tau),p(\tau)\rangle= (q(\tau)+\langle \hat{w}(\tau),p(\tau)\rangle)\langle\nabla_u\hat{L}(\tau),\hat{w}(\tau)\rangle.
\end{equation}
Since $q(\tau)\leq 0$, we have $\langle \hat{w}(\tau),p(\tau)\rangle>0$. From (\ref{P44}) we get
$$
 \langle\nabla_u\hat{L}(\tau),\hat{w}(\tau)\rangle=\frac{\langle \hat{w}(\tau),p(\tau)\rangle}{q(\tau)+\langle \hat{w}(\tau),p(\tau)\rangle}( \hat{L}(\tau)+\beta)\geq  \hat{L}(\tau)+\beta.
$$
From this we obtain
$$
\beta\leq \langle\nabla_u\hat{L}(\tau),\hat{w}(\tau)\rangle-\hat{L}(\tau)\leq \sigma (|\hat{w}(\tau)|).
$$
Hence, we have
$$
|\hat{w}(\tau)|\geq\sigma^{-1}(\beta)= \sigma^{-1}\left(\sigma\left(\frac{c+1}{T_0}\right)\right)=\frac{c+1}{T_0}.\;\;\;\qed
$$

\begin{lemma}
\label{lem4}
If $q(\tau_1)= 0$ and $q(\tau)< 0$, $\tau\in ]\tau_1,\tau_2]$, then $|q(\tau)|/|p(\tau)|\leq \bar{\gamma}$, $\tau\in [\tau_1,\tau_2]$.
\end{lemma}

{\it Proof}. Since $q(\tau)< 0$, $\tau\in ]\tau_1,\tau_2]$, we have
$$
\frac{q(\tau)+\langle p(\tau),\hat{w}(\tau)\rangle}{|p(\tau)|( \hat{L}(\tau)+\beta)}\leq \frac{|\hat{w}(\tau)|}{\theta(|\hat{w}(\tau)|)+\beta}\leq\eta.
$$
From (\ref{P1}), (\ref{P2}), and condition (C3) we get
\begin{eqnarray}
&& \left|\frac{dq(\tau)}{d\tau}\right| \leq \frac{ (q(\tau)+\langle p(\tau),\hat{w}(\tau)\rangle)|\hat{L}_t(\tau)|}{( \hat{L}(\tau)+\beta)^2}\nonumber\\
&&\leq\eta |p(\tau)| \frac{ |\hat{L}_t(\tau)|}{ \hat{L}(\tau)+\beta}\leq\eta\xi|p(\tau)|,\label{u1}\\
&& \left|\frac{dp(\tau)}{d\tau}\right| \leq \frac{ (q(\tau)+\langle p(\tau),\hat{w}(\tau)\rangle)|\nabla_x\hat{L}(\tau)|}{( \hat{L}(\tau)+\beta)^2}\nonumber\\
&&\leq\eta |p(\tau)| \frac{ |\nabla_x\hat{L}(\tau)|}{ \hat{L}(\tau)+\beta}\leq\eta\xi|p(\tau)|.\label{u2}
\end{eqnarray}
whenever $\tau\in [\tau_1,\tau_2]$. 
From this we obtain
$$
\frac{d}{d\tau}\frac{|q(\tau)|}{|p(\tau)|}\leq\frac{|dq(\tau)/d\tau||p(\tau)|+|q(\tau)||dp(\tau)/d\tau|}{|p(\tau)|^2}
\leq\eta\xi\left( 1+\frac{|q(\tau)|}{|p(\tau)|}\right).
$$
Since $q(\tau_1)=0$, applying the Gronwall inequality we have
\begin{equation}
\label{utt}
\frac{|q(\tau)|}{|p(\tau)|}\leq e^{\eta\xi (\tau_2-\tau_1)}-1.
\end{equation}
Observe that
$$
\tau_2-\tau_1=\int_{\hat{t}(\tau_1)}^{\hat{t}(\tau_2)}( L(t,\hat{x}(t),\dot{\hat{x}}(t))+\beta)dt\leq  c+(\hat{t}(\tau_2)-\hat{t}(\tau_1))\beta, 
$$
and $\hat{t}(\tau_2)-\hat{t}(\tau_1)\leq T_0$. Indeed, if $q(\tau)\leq 0$ on a time interval $[\tau_1,\tau_2]$ and $t\in [\hat{t}(\tau_1),\hat{t}(\tau_2)]$, then by Lemma \ref{lem2} we have
$$
|\dot{\hat{x}}(t)|=|\hat{w}(\hat{\tau}(t))|\geq \frac{c+1}{T_0}.
$$ 
The inequality $\hat{t}(\tau_2)-\hat{t}(\tau_1)> T_0$ contradicts (\ref{ner1}). Thus $\tau_2-\tau_1\leq  c+T_0\beta$. Combining this with (\ref{utt}) and (\ref{gam}) we obtain the result. \qed

\vspace{5mm}

{\it End of the proof of Theorem \ref{th1}}. Let us consider $\tau\in [0,\hat{T}]$ such that $q(\tau)\geq 0$.
From condition (C2) we have
$$
\hat{L}(\tau)-\langle \nabla_u\hat{L}(\tau),\hat{w}(\tau)\rangle +\frac{\mu}{2}|\hat{w}(\tau)|^2\leq \Lambda_0.
$$
From this and (\ref{P44}) we obtain
$$
\frac{\mu}{2}|\hat{w}(\tau)|^2\leq \Lambda_0-\hat{L}(\tau)+\langle \nabla_u\hat{L}(\tau),\hat{w}(\tau)\rangle=
$$
$$
\Lambda_0-\hat{L}(\tau)+\frac{\langle p(\tau),\hat{w}(\tau)\rangle}{q(\tau)+\langle p(\tau),\hat{w}(\tau)\rangle}(\hat{L}(\tau)+\beta).
$$
If $\langle p(\tau),\hat{w}(\tau)\rangle> 0$, then we have
$$
\frac{\mu}{2}|\hat{w}(\tau)|^2\leq \Lambda_0-\hat{L}(\tau)+\hat{L}(\tau)+\beta=\Lambda_0+\beta.
$$
Hence
\begin{equation}
\label{nn1}
|\hat{w}(\tau)|\leq\sqrt{\frac{2}{\mu}\left( \Lambda_0+\beta\right)}.
\end{equation}
If $\langle p(\tau),\hat{w}(\tau)\rangle\leq 0$, then from Lemma \ref{lem00} we get
$$
\frac{\mu}{2}|\hat{w}(\tau)|^2\leq \Lambda_0-\hat{L}(\tau)\leq \Lambda_0+\Lambda_1 |\hat{w}(\tau)| -\frac{\mu}{2} |\hat{w}(\tau)|^2.
$$
Thus we obtain 
\begin{equation}
\label{nn2}
|\hat{w}(\tau)|\leq \frac{\Lambda_1+\sqrt{\Lambda_1^2+4\mu\Lambda_0}}{2}.
\end{equation}

Now, let us consider $\tau\in [0,\hat{T}]$ such that $q(\tau)< 0$. By Lemma \ref{lem4} $q(\tau)>-\bar{\gamma} |p(\tau)|$. From the maximum principle we have
$$
\frac{q(\tau)+\langle p(\tau),\hat{w}(\tau)\rangle}{ \hat{L}(\tau)+\beta}\geq
\frac{q(\tau)+(\bar{\gamma}+1)\langle p(\tau),p(\tau)/|p(\tau)|\rangle}{ {L}(\hat{t}(\tau),\hat{y}(\tau),(\bar{\gamma}+1)p(\tau)/|p(\tau)|)+\beta}
$$
\begin{equation}
\label{1000}
\geq
\frac{|p(\tau)|}{ {L}(\hat{t}(\tau),\hat{y}(\tau),(\bar{\gamma}+1)p(\tau)/|p(\tau)|)+\beta}
\geq\frac{|p(\tau)|}{\Lambda_2+\beta}.
\end{equation}
Recall that $\varrho>0$ is a number such that for all $r>\varrho$ (\ref{rh}) is satisfied.  If $|w|>\varrho$, then by Lemma \ref{lem00}  we get
$$
\frac{q(\tau)+\langle p(\tau),w\rangle}{ {L}(\hat{t}(\tau),\hat{y}(\tau),w)+\beta}
\leq \frac{|p(\tau)||w|}{(\frac{\mu}{2}|w|^2-\Lambda_1|w|-\Lambda_0)+\beta}\leq\frac{|p(\tau)|}{\Lambda_2+\beta}.
$$
From  (\ref{1000}) we see that $w\neq \hat{w}(\tau)$. Thus $|\hat{w}(\tau)|\leq \varrho$. Combining this with (\ref{nn1}) and (\ref{nn2}) we obtain the result. \qed

\subsection{Proof of Theorem \ref{th3}}

 Let us consider the following problem
\begin{eqnarray}
\label{30}
&& \int_0^1 L(t,x(t),\dot{x}(t))dt \rightarrow\inf, \nonumber\\
&& \dot{x}\in \mathcal{U}_N, \\
&& x(0)=0,\;\;\; Ax(1)\leq b. \nonumber
\end{eqnarray}
Obviously this problem has a solution $\tilde{x}(\cdot)$. Let $\bar{x}(\cdot)\in AC([0,1],R^n)$ be the function  defined by
\begin{eqnarray}
&& \bar{x}(0)=0,\nonumber\\
&& \bar{x}(t)=\bar{x}(k\tau)+\frac{t-k\tau}{\tau}\int_{k\tau}^{(k+1)\tau}\dot{\hat{x}}(s)ds,\\
&& t\in [k\tau, (k+1)\tau), \ \ \ k\in\overline{1,N-1}, \ \ \ \tau=1/N.\nonumber
\end{eqnarray}
It satisfies the conditions $\bar{x}(k\tau )=\hat{x}(k\tau )$, $k=\overline{0,N}$, and $|\bar{x}(t)-\hat{x}(t)|\leq \ell\tau$, $t\in [0,1]$.
Since the function $ L(t,x,\cdot)$ is convex and continuous, and the function $L(\cdot,\cdot, u)$ is Lipschitzian with the constant $K_L$ in the set $\{ (t,x,u)\mid t\in [0,1], |x|\leq\ell,  |u|\leq\ell\}$, we have
$$
\int_0^1 L(t,\bar{x}(t),\dot{\bar{x}}(t))dt
=
\sum_{k=0}^{N-1}
\int_{k\tau}^{(k+1)\tau}L\left(t,\bar{x}(t),\frac{1}{\tau}\int_{k\tau}^{(k+1)\tau}\dot{\hat{x}}(s)ds\right) dt
$$
$$
\leq
\sum_{k=0}^{N-1}
\int_{k\tau}^{(k+1)\tau}L\left(k\tau,\hat{x}(k\tau),\frac{1}{\tau}\int_{k\tau}^{(k+1)\tau}\dot{\hat{x}}(s)ds\right) dt+\tau K_L(\ell +1)
$$
$$
\leq
\sum_{k=0}^{N-1}
\int_{k\tau}^{(k+1)\tau}\frac{1}{\tau}\int_{k\tau}^{(k+1)\tau}L(k\tau,\hat{x}(k\tau),\dot{\hat{x}}(s))dsdt+\tau K_L(\ell +1)
$$
$$
=\sum_{k=0}^{N-1}
\int_{k\tau}^{(k+1)\tau}L(k\tau,\hat{x}(k\tau),\dot{\hat{x}}(t))dt+\tau K_L(\ell +1)
$$
$$
\leq\sum_{k=0}^{N-1}
\int_{k\tau}^{(k+1)\tau}L(t,\hat{x}(t),\dot{\hat{x}}(t))dt+2\tau K_L(\ell +1)
$$
$$
=
\int_0^1L(t,\hat{x}(t),\dot{\hat{x}}(t))dt+2\tau K_L(\ell +1).
$$
Since
$$
\int_0^1 L(t,\hat{x}(t),\dot{\hat{x}}(t))dt
\leq\int_0^1 L(t,\tilde{x}(t),\dot{\tilde{x}}(t))dt
\leq\int_0^1 L(t,\bar{x}(t),\dot{\bar{x}}(t))dt,
$$
we obtain
$$
\left|
\int_0^1 L(t,\tilde{x}(t),\dot{\tilde{x}}(t))dt-
\int_0^1 L(t,\hat{x}(t),\dot{\hat{x}}(t))dt
\right|\leq 2\tau K_L(\ell +1).
$$
From this and the inequality
$$
N>\frac{4(1+\ell)K_L}{\varepsilon},
$$
we have
$$
\left|
\int_0^1 L(t,\tilde{x}(t),\dot{\tilde{x}}(t))dt-
\int_0^1 L(t,\hat{x}(t),\dot{\hat{x}}(t))dt
\right|\leq \frac{\varepsilon}{2}.
$$
Using the path-following method we can find $\breve{u}\in \mathcal{U}_N$ such that
$$
\breve{x}(t)=\int_0^t\breve{u}(s)ds.
$$
is an admissible solution to problem (\ref{30}) satisfying 
$$
 \left|
\int_0^1 L(t,\tilde{x}(t),\dot{\tilde{x}}(t))dt-
\int_0^1 L\left(t,\breve{x}(s),\breve{u}(t)\right)dt
\right|\leq \frac{\varepsilon}{2}.
$$
Thus, we get 
$$
\left|
\int_0^1 L\left(t,\breve{x}(s),\breve{u}(t)\right)dt-
\int_0^1 L(t,\hat{x}(t),\dot{\hat{x}}(t))dt
\right|\leq \varepsilon.
$$

Since $|\dot{\hat{x}}|\leq\ell$ we can use the necessary conditions of optimality:
\begin{eqnarray}
&& \frac{d}{dt}\nabla_u L(t,\hat{x}(t),\dot{\hat{x}}(t))=\nabla_xL(t,\hat{x}(t),\dot{\hat{x}}(t)),\label{nc1}\\
&& \langle\nabla_u L(1,\hat{x}(1),\dot{\hat{x}}(1)),z-\hat{x}(1)\rangle\geq 0,\;\; z\in S.\label{nc2}
\end{eqnarray}
From (\ref{sc}) we obtain
$$
\int_0^1 L(t,\hat{x}(t),\dot{\hat{x}}(t))dt
$$
$$
+\int_0^1 (\langle\nabla_x L(t,\hat{x}(t),\dot{\hat{x}}(t)), \breve{x}(t)-\hat{x}(t)\rangle
+\langle\nabla_u L(t,\hat{x}(t),\dot{\hat{x}}(t)), \dot{\breve{x}}(t)-\dot{\hat{x}}(t)\rangle )dt
$$
$$
+\int_0^1 \frac{\mu}{2}(|\breve{x}(t)-\hat{x}(t)|^2+|\dot{\breve{x}}(t)-\dot{\hat{x}}(t)|^2)dt\leq
\int_0^1 L(t,\breve{x}(t),\dot{\breve{x}}(t))dt.
$$
Integrating the third term in the left-hand side of the inequality by parts and using (\ref{nc1}) and (\ref{nc2}), we get
$$
\int_0^1 (|\breve{x}(t)-\hat{x}(t)|^2+|\dot{\breve{x}}(t)-\dot{\hat{x}}(t)|^2)dt
$$
$$
\leq
\frac{2}{\mu}
\left|
\int_0^1 L(t,\breve{x}(t),\dot{\breve{x}}(t))dt-
\int_0^1 L(t,\hat{x}(t),\dot{\hat{x}}(t))dt
\right|\leq \frac{2}{\mu}\varepsilon.
$$

By Proposition \ref{proposition5}, the number of iterations needed to find $\breve{u}(\cdot)$ does not exceed
$$
\frac{\ln\left( 2\frac{(1+\kappa)(m+nN+1)+(\kappa+\sqrt{m+nN+1})\kappa}{\gamma(1-2\kappa)\varepsilon}\|v\|_{(\dot{\tilde{x}}, \tilde \sigma)}^{\mathcal{F}}\right)}{\ln\left (1+\frac{\gamma}{\kappa+\sqrt{m+nN+1}}\right)}+1 .
$$
The term $\|v\|_{(\dot{\tilde{x}}, \tilde \sigma)}^{\mathcal{F}}$ can be estimated by $\frac{1}{\sqrt{2}}\bar{\sigma}$ (see (\ref{sig})). This ends the proof.

\section{Conclusion}

In this work we analyzed a convex  problem of calculus of variations with polyhedral end-point constraints.  Our objective was to get complexity bounds for a path-following method applied to the problem. To this end we deduced explicit estimate for the Lipschitz constant of solution to the problem of calculus of variations. This allowed us to compute the time interval partition diameter needed to construct piecewise linear approximation of solution with given accuracy and to reduce the original problem to a convex programming one.  Then we showed that the path-following method finds an admissible piecewise linear solution approximating the solution to the original problem with given accuracy. The estimate for the  number of iterations obtained in the paper depends only on the data of the original problem of calculus of variations.

\subsection*{Acknowledgements}
The authors are grateful to Delfim Torres for bibliographical support. This work was partially supported by project PTDC/EEI-AUT/2933/2014 (TOCCATA), funded by Project 3599 - Promover a 
Produ\c c\~ao Cient\'\i fica e Desenvolvimento Tecnol\'ogico
e a Constitui\c c\~ao de Redes Tem\'aticas (3599-PPCDT) and FEDER funds through
COMPETE 2020, Programa Operacional Competitividade e Internacionaliza\c c\~ao (POCI),
and by national  funds through Funda\c c\~ao para a Ci\^encia e a Tecnologia (FCT). The work of Miguel Oliveira was supported by FCT through the PhD fellowship SFRH/BD/111854/2015.

\end{document}